\documentclass[11pt,leqno]{amsart}
\usepackage{amsmath,amsfonts,amssymb,amscd,amsthm,amsbsy,upref}
\newtheorem{thm}{Theorem}
\newtheorem{lem}[thm]{Lemma}
\newtheorem*{sublem}{Sublemma}

\newtheorem{prop}[thm]{Proposition}
\newtheorem*{mainthm}{Main Theorem}
\theoremstyle{definition}
\newtheorem{problem}[thm]{Problem}
\def\ds{\displaystyle}

\def\nat{{\mathbb N}}
\def\real{{\mathbb R}}
\def\A{{\mathcal A}}
\def\F{{\mathcal F}}
\def\S{{\mathcal S}}
\def\ep{\varepsilon}
\def\w{\omega}
\begin{document}
\author{W.~B.~Johnson \and E.~Odell}
\title{The Diameter of the Isomorphism Class of a Banach Space}
\address{Department of Mathematics, Texas A\&M University, 
College Station, TX 77843-3368}
\email{johnson@math.tamu.edu}
\address{Department of Mathematics, The University of Texas at Austin,
1 University Station C1200,
Austin, TX 78712-0257}
\email{odell@math.utexas.edu}

\begin{abstract}
We prove that if $X$ is a separable infinite dimensional Banach space then 
its isomorphism class has infinite diameter with respect to the 
Banach-Mazur distance. 
One step in the proof is to show that if $X$ is elastic then $X$ contains 
an isomorph of $c_0$. 
We call $X$ elastic if for some $K <\infty$ for every Banach space $Y$
which embeds into $X$,  the space $Y$ is $K$-isomorphic to a subspace of $X$.
We also prove that if $X$ is a separable Banach space such that for some 
$K<\infty$ every isomorph of $X$ is $K$-elastic then $X$ is finite dimensional.
\end{abstract}
\thanks{Johnson was supported
in part by NSF DMS-0200690,  Texas Advanced Research 
Program 010366-163, and the U.S.-Israel
Binational Science Foundation. Odell was supported
in part by NSF DMS-0099366 
and was a participant in the NSF supported Workshop in
Linear Analysis and Probability at Texas A\&M University.\hfill\break MR
subject classification: 46G05, 46T20.}
\dedicatory{September 23, 2003}

\maketitle
\baselineskip=18pt               

\begin{section}{Introduction.}\label{intro}


Given a Banach space $X$, let $D(X)$ be the diameter in the
Banach-Mazur distance of the class of all Banach spaces which are
isomorphic to $X$; that is, 
\[ 
D(X) = \sup \{d(X_1,X_2) : X_1, X_2 \mathrm{\ are\ isomorphic\
to }
\, X\}
\]
where  $d(X_1,X_2) $ is the infimum over all  isomorphisms
$T$ from $X_1$ onto $X_2$ of $\|T\|\cdot \|T^{-1}\|$.  It is well known that if
$X$ is finite (say, $N$) dimensional, then $c N\le D(X) \le N$ for some
positive constant $c$ which is independent of $N$.  The upper bound is an 
immediate
consequence of the classical result (see e.g. \cite[p.~54]{tomczakbook}) that
$d(Y,\ell_2^N)\le
\sqrt N$ for every $N$ dimensional space $Y$.  The lower bound is due to Gluskin
\cite{gluskin}, \cite[p.~283]{tomczakbook}.  

It is natural to conjecture that $D(X) $ must be infinite when $X$ is infinite
dimensional, but this problem remains open. As far as we know, this problem
was first raised in print in the 1976 book of  J.~J.~Sch\"affer 
\cite[p. 99]{S}.   The problem was recently brought to  the attention of the
authors by V.~I.~Gurarii, who checked that every infinite dimensional
super-reflexive space as well as each of the common classical Banach spaces
has an isomorphism class  whose diameter is infinite. To see these cases, note
that if
$X$ is infinite dimensional and $E$ is any finite dimensional space, then it 
is clear
that $X$ is isomorphic to $E\oplus_2 X_n$ for some space $X_n$.  Therefore, if
$D(X)$ is finite, then $X$ is finitely complementably universal; that is,
there is a constant $C$ so that every finite dimensional space is $C$-isomorphic
to a $C$-complemented subspace of $X$.  This implies that $X$ cannot have non
trivial type or non trivial cotype or local unconditional structure or numerous
other structures.  In particular,
$X$ cannot be any of the classical spaces or be super-reflexive.

 In his unpublished 1968 thesis \cite{Mc}, McGuigan 
conjectured that $D(X)$ must  be larger than one when dim\,$X>1$.  
Sch\"affer \cite[p. 99]{S} derived that $D(X)\ge 6/5$ when dim\,$X>1$ as a
consequence of other geometrical results contained in \cite{S}, but one can  prove
directly that $D(X)\ge \sqrt{2}$.  Indeed, it  is clearly enough to get an
appropriate lower bound on the  Banach-Mazur distance between
$X_1:= Y \oplus_1 \real$ and 
$X_2:= Y \oplus_2 \real$    when $Y$ is a non zero Banach space.  Now $X_1$ has
a one dimensional subspace for which every two dimensional superspace is
isometric to $\ell_1^2$.  On the other hand, every one dimensional subspace of
$X_2$ is contained in a two dimensional superspace which is isometric to
$\ell_2^2$.  It follows that $d(X_1,X_2)\ge d(\ell_1^2,\ell_2^2)=\sqrt{2}$.

The Main Theorem in this paper is a solution to Sch\"affer's problem for
separable Banach spaces:

\begin{mainthm}\label{thm:main}
If $X$ is a separable infinite dimensional Banach space, then $D(X)=\infty$.
\end{mainthm}


Part of the work for proving the Main Theorem involves showing that 
if $X$ is
separable and $D(X)<\infty$, then $X$ contains an isomorph of $c_0$. 
This proof is inherently non local in nature, and, strangely enough,  
local considerations, such as those mentioned earlier which yield 
partial results, play no role in our proof.
We do not see how to prove that a non separable space $X$ for 
which  $D(X)<\infty$ must contain an isomorph of $c_0$.  
Our proof requires Bourgain's index theory which in turn requires 
separability.

Our method of proof involves the concept of an {\sl elastic Banach space.} 
Say that $X$ is {\em $K$-elastic\/} provided that if a Banach space 
$Y$ embeds into $X$ then $Y$ must {\em $K$-embed into\/} $X$ 
(that is, there is an isomorphism $T$  from $Y$ into $X$ with 
$$\|y\| \le \|Ty\| \le K\|y\|$$ 
for all $y\in Y$).
This is the same (by Lemma~\ref{lem1}) as saying that every space 
isomorphic to $X$ must $K$-embed  into $X$. 
$X$ is said to be {\em elastic\/} if it is $K$-elastic for some $K<\infty$.

Obviously, if $D(X)<\infty$
then $X$ as well as every isomorph of $X$ is   $D(X)$-elastic. 
Thus the Main Theorem
is an immediate consequence of

\begin{thm}\label{thm:elastic} If $X$ is a separable Banach space and there is a $K$
so that every isomorph of $X$ is $K$-elastic, then $X$ is finite dimensional.
\end{thm}
 
A key step in our argument involves showing that an elastic space $X$ admits 
a normalized weakly null sequence having a spreading model not equivalent 
to either the unit vector basis of $c_0$ or $\ell_1$.   To achieve this we 
first prove (Theorem~\ref{thm6})
that if $X$ is elastic then 
$c_0$ embeds into $X$. It is reasonable to conjecture that an elastic 
infinite dimensional separable Banach space must contain an isomorph of $C[0,1]$. 
Theorem \ref{thm:elastic} would be an immediate consequence of this conjecture and
the ``arbitrary distortability" of $C[0,1]$ proved in \cite{LP}.  Our derivation of
Theorem \ref{thm:elastic} from Theorem~\ref{thm6} uses ideas from \cite{LP} as well
as  \cite{MR}.

With the letters $X,Y,Z,\ldots$ we will denote separable infinite dimensional real 
Banach spaces unless otherwise indicated. 
$Y\subseteq X$ will mean that $Y$ is a closed (infinite dimensional) 
subspace of $X$. The closed linear span of the set $A$ is denoted $[A]$.  We
use standard Banach space theory terminology, as can be found in \cite{LT}. 
 The 
material we use on spreading models can be found  in
\cite{BL}. For simplicity we assume real scalars, but all proofs can easily be
adapted for complex Banach spaces.
\end{section}

\begin{section}{The Main Result}\label{result}

The following well known  elementary lemma shows that the two definitions of
elastic mentioned in Section \ref{intro} are equivalent.

\begin{lem}\label{lem1} 
Let $Y\subseteq (X,\|\cdot\|)$ and let $|\cdot|$ be an equivalent norm 
on $(Y,\|\cdot\|)$. 
Then $|\cdot|$ can be extended to an equivalent norm on $X$.
\end{lem}

\begin{proof}
There exist positive reals $C$ and $d$ with $d\|y\| \le |y| \le C\|y\|$ 
for $y\in Y$. 
Let $F\subseteq CB_{X^*}$ be a set of Hahn-Banach extensions of all elements of 
$S_{(Y^*,|\cdot|)}$ to all of $X$. 
For $x\in X$ define 
$$|x| = \sup \big\{ |f(x)| : f\in F\big\} \vee d\|x\|\ .
\qquad \qed$$
\renewcommand{\qed}{}
\end{proof}

Let $n\in \nat$ and $K<\infty$.
We shall call a basic sequence $(x_i)$ {\em block $n$-unconditional with 
constant $K$} if every block basis $(y_i)_{i=1}^n$ of $(x_i)$ is 
$K$-unconditional; that is, 
\[
\|\sum_{i=1}^n \pm a_i y_i\|\le K \|\sum_{i=1}^n   a_i y_i\|
\]
for all scalars $(a_i)_{i=1}^n$ and all choices of $\pm$.

The next lemma is essentially contained in \cite{LP}.  In fact, by using the
slightly more involved  argument in \cite{LP}, the conclusion
``with constant $2$"  can be changed to ``with constant
$1+\ep$", which implies that the constant in the conclusion of Lemma
\ref{lem2} can be changed from $2+\ep$ to $1+\ep$. 

\begin{lem}\label{blockUncon}
Let $X$ be a Banach space with a basis $(x_i)$. For every $n$ there is
an equivalent norm $|\cdot|_n$ on $X$ so that in $(X,|\cdot|_n)$, $(x_i)$
is block $n$-unconditional with constant $2$.
\end{lem}
\begin{proof}
Let $(P_n)$ be the sequence of basis projections associated with $(x_n)$. 
We may assume, by passing to an equivalent norm on $X$, that $(x_n)$ is 
bimonotone and hence $\|P_j-P_i\|=1$ for all $i<j$.
Let $\S_n$ be the class of operators $S$ on $X$ of the form 
$S= \sum_{k=1}^m (-1)^k (P_{n_k} - P_{n_{k-1}})$ where 
$0\le n_0 < \cdots < n_m$ and $m\le n$. 
Define 
$$|x|_n
 := \sup \{ \|Sx\| : S\in \S_n\}\ .$$
Thus $\|x\| \le |x|_n \le n\|x\|$ for $x\in X$. 
It suffices to show that for $S\in \S_n$, $|S|_n \le 2$. 
Let $x\in X$ and $|Sx|_n= \|TSx\|$ for some $T\in S_n$.
Then since $TS\in \S_{2n} \subseteq \S_n + \S_n$,
$$\|TSx\| \le 2|x|_n\, .$$
\end{proof} 

\begin{lem}\label{lem2} 
For every separable Banach space $X$, $n\in\nat$, and $\ep>0$, there
exists an equivalent norm
$|\cdot|$  on $X$ so that every normalized weakly null sequence in $X$ admits
a  block $n$-unconditional subsequence with constant $2+\ep$.
\end{lem}

\begin{proof} 
Since $C[0,1]$ has a basis, the lemma follows from Lemma \ref{blockUncon} 
and the
the classical fact that every separable Banach space  
$1$-embeds 
 into $C[0,1]$.  
\end{proof}

Lemma \ref{lem2} is false for some non separable spaces.  
Partington \cite{P} and Talagrand \cite{T}   proved
that every isomorph of $\ell_\infty$  contains, for every $\ep>0$, a 
$1+\ep$-isometric copy of
$\ell_\infty$ and hence of every separable Banach space.

Our next lemma is an extension of the Maurey-Rosenthal construction 
\cite{MR}, or rather the footnote to it given by one of the authors (Example
3 in  \cite{MR}).   We first recall the construction of {\em spreading
models}. If $(y_n)$ is a normalized  basic sequence then, given $\ep_n
\downarrow 0$, one can use Ramsey's theorem and a diagonal argument to find
a subsequence
$(x_n)$ of $(y_n)$ with the following property.  For all $m$ in $\nat$ and
$(a_i)_{i=1}^m \subset [-1,1]$, if $m\le i_1<\dots<i_m$ and $m\le
j_1<\dots<i_m$, then
$$
\Bigg{|} \Big{\|}\sum_{k=1}^m a_k x_{i_k}\Big{\|} - \Big{\|}\sum_{k=1}^m a_k
x_{j_k}\Big{\|}\Bigg{|} < \ep_m.
$$
It follows that for all $m$ and $(a_i)_{i=1}^m \subset \real$,
$$
\lim_{i_1\to\infty} \dots\lim_{i_m\to\infty} 
\Big{\|}\sum_{k=1}^m a_k x_{i_k}\Big{\|} \equiv \Big{\|}\sum_{k=1}^m a_k
\tilde x_k\Big{\|}
$$
exists.  
The sequence $(\tilde x_i)$ is then a basis for the completion of
$(\text{span}\, (\tilde{x}_{i }), \|\cdot\|)$ and $(\tilde{x}_{i})$ is
called a  {\em spreading model} of $(x_i)$.
  If $(x_i)$ is weakly null, then $(\tilde{x}_i)$ is $2$-unconditional. 
One shows this by checking that $(\tilde{x}_i)$ is suppression
$1$-unconditional, which means that for all scalars $(a_i)_{i=1}^m$ and
$F\subset \{1,\dots, m\}$,
$$
\|\sum_{i\in F} a_i \tilde{x}_i\|\le \|\sum_{i=1}^m a_i \tilde{x}_i\|.
$$
Also, $(x_i)$ is $1$-subsymmetric, which means that for all scalars
$(a_i)_{i=1}^m$ and all $n(1)<\dots n(m)$,
$$
\|\sum_{i=1}^m a_i \tilde{x}_i\| \le
\|\sum_{i=1}^m a_i \tilde{x}_{n(i)}\|.
$$
It is not difficult to see that, when $(x_i)$ is weakly null,
$(\tilde{x}_i)$ is not equivalent to the unit vector basis of $c_0$
(respectively, $\ell_1$) if and only if
$\lim_m \|\sum_{i=1}^m  \tilde{x}_i\| =\infty$ (respectively,
$\lim_m \|\sum_{i=1}^m  \tilde{x}_i\|/m = 0$).  
All of these facts can be found in \cite{BL}. 

\begin{lem}\label{lem3}
Let $(x_n)$ be a normalized weakly null basic sequence with spreading 
model $(\tilde x_n)$. 
Assume that $(\tilde x_n)$ is not equivalent to either the unit vector 
basis of $\ell_1$ or the unit vector basis of $c_0$. 
Then for all $C<\infty$ there exist $n\in\nat$, 
a subsequence $(y_i)$ of $(x_i)$, and  
an equivalent norm $|\cdot|$ on $[(y_i)]$ so that $(y_i)$ is  
$|\cdot|$-normalized 
and no subsequence of $(y_i)$ is block $n$-unconditional with constant 
$C$ for the norm $|\cdot|$.
\end{lem}

\begin{proof}
Recall that if $(e_i)_1^n$ is normalized and $1$-subsymmetric then 
$\|\sum_1^n e_i\|\, \|\sum_1^n e_i^*\| \le 2n$ where $(e_i^*)_1^n$ 
is biorthogonal to $(e_i)_1^n$ \cite[p.118]{LT}. 
Thus $e = \frac{\sum_1^n e_i}{\|\sum_1^n e_i\|}$ is normed by 
$f= \frac{\|\sum_1^n e_i\|}{n} \sum_{i=1}^n e_i^*$, 
precisely $f(e) =1 = \|e\|$, and $\|f\| \le 2$. 
These facts allow us to deduce that there is a subsequence $(y_i)$ 
of $(x_i)$ so that if 
$F\subseteq \nat$ is {\em admissible} (that is, $|F| \le \min F $) then 
$$f_F \equiv \frac{\|\sum_{i\in F} y_i\|}{|F|} 
\Big( \sum_{i\in F} y_i^*\Big)$$
satisfies  $\|f_F\| \le 5$ and $f_F (y_F)=1$, where 
$$y_F \equiv \frac{\sum_{i\in F} y_i}{\|\sum_{i\in F} y_i\|}\ .$$

Indeed, $(\tilde x_i)$ is $1$-subsymmetric and suppression 
$1$-unconditional
(since $(x_i)$ is weakly null).  
Given $1/2>\ep >0$ we can  find $(y_i)\subseteq (x_i)$ so that if 
$F\subseteq \nat$ is admissible then $(y_i)_{i\in F}$ is 
$1+\ep$-equivalent to $(\tilde x_i)_{i=1}^{|F|}$. 
Furthermore we can choose $(y_i)$ so that if $F$ is admissible
then for $y=\sum a_iy_i$, 
$\|\sum_{i\in F} a_i y_i\| \le (2+\ep) \|y\|$ 
(\cite{MR}; for a proof see \cite{O} or \cite{BL}).
Hence $\|f_F\| \le (2+\ep)\|f_F\mid_{[ y_i]_{i\in F}}\| < 5$
for sufficiently small $\ep$ by our above remarks.

We are ready to produce a Maurey-Rosenthal type renorming. 
Choose $n$ so that $n>7C$ and let $\ep>0$ satisfy $n^2\ep <1$. 
We choose a subsequence $M= (m_j)_{j=1}^\infty$ of $\nat$ so that 
$m_1=1$ and for $i\ne j$ and for all admissible sets $F$ and $G$ with 
$|F|=m_i$ and $|G| = m_j$, 
\begin{itemize}
\item[a)] $\ds  \frac{\|\sum_{i\in F}y_i\|}{\|\sum_{i\in G} y_i\|}
<\ep$ , if  $m_i <m_j$ and 
\vskip10pt
\item[b)] $\ds \frac{\|\sum_{i\in F} y_i\|}{\|\sum_{i\in G}y_i\|} \ 
\frac{m_j}{m_i} < \ep$ , if $m_i >m_j\ .$
\end{itemize}

Indeed, we have chosen $(y_i)$ so that 
$$\frac12 \Big\| \sum_{i=1}^{|F|} \tilde x_i\Big\| 
\le \Big\| \sum_{i\in F} y_i\Big\| 
\le 2\Big\| \sum_{i=1}^{|F|} \tilde x_i\Big\|$$
and similarly for $G$. 
Since $(\tilde x_i)$ is not equivalent to the unit vector basis of $c_0$ 
(and is unconditional) 
$\lim_m \| \sum_1^m \tilde x_i\| = \infty$
so that a) will be satisfied if
 $(m_k)$  increases sufficiently rapidly. 
Furthermore, since $(\tilde y_i)$ is  
not equivalent to the unit  vector basis of $\ell_1$,   $\lim_m
\frac{\|\sum_1^m \tilde y_i\|}{m} = 0$  and so  b) can also be achieved.

For $i\in\nat$ set $\A_i = \{y_F : F$ is admissible and $|F|=m_i\}$ and 
$\A_i^* = \{f_F:F$ is admissible and $|F|=m_i\}$. 
Let $\phi$ be an injection into $M$ from the collection of all 
$(F_1,\ldots,F_i)$ where  $i<n$ and $F_1 < F_2 <\cdots < F_i$ are finite 
subsets of $\nat$. 
Let 
$$\F = \bigg\{ \sum_{i=1}^n f_{F_i} : F_1 < \cdots < F_n,\ 
|F_1| = m_1 =1,\ |F_{i+1}| =  {\phi (F_1,\ldots,F_i)} \text{ for } 
1\le i<n\bigg\}$$
For $y\in [(y_i)]$ let 
$$\|y\|_{\F} = \sup \Big\{ |f(y)| : f\in \F\Big\}$$ 
and set 
$$|y| = \|y\|_{\F} \vee \ep \|y\|\ .$$
This is an equivalent norm since for $f\in\F$, $\|f\| \le 5n$.

Note that if $f_F \in \A_i^*$ and $y_G \in \A_j$ with $i\ne j$ then 
$$|f_F(y_G)| = \frac{\|\sum_{i\in F} y_i\|}{m_i} 
\sum_{i\in F} y_i^* \left( \frac{\sum_{i\in G} y_i}{\|\sum_{i\in G} y_i\|} 
\right)
\le \frac{\|\sum_{i\in F} y_i\|}{\|\sum_{i\in G} y_i\|} 
\frac{m_i\wedge m_j}{m_i}\ .$$
If $m_i < m_j$ then $|f_F (y_G)|<\ep$ by a). 
If $m_i > m_j$ then $|f_F(y_G)| <\ep$ by b). 

It follows that if $y= \sum_{i=1}^n y_{F_i}$ and 
$f= \sum_{i=1}^n f_{F_i} \in \F$ then $|y| \ge f(y)=n$ and if 
$z= \sum_{i=1}^n (-1)^i y_{F_i}$, then for all 
$g = \sum_{j=1}^n f_{G_j} \in \F$, $|g(z)| \le 6 + n^2 \ep <7$. 
Indeed, we may assume that $g\not= f$ and if $F_1\ne G_1$ then $|G_j| \ne
|F_i|$ for all
$1<i\le n$ and 
$1\le j\le n$ and so by a), b), 
$$|g(z)| = \bigg| \sum_{j=1}^n f_{G_j} \bigg( \sum_{i=1}^n (-1)^i y_{F_i}
\bigg)\bigg| 
\le \sum_{j=1}^n \sum_{i=1}^n |f_{G_j} (y_{F_i})| < n^2 \ep\ .$$
Otherwise there exists $1\le j_0 <n$ so that $F_j = G_j$ for $j\le j_0$, 
$|F_{j_0+1}| = |G_{j_0+1}|$ and $|F_i| \ne |G_j|$ for 
$j_0 + 1 < i,j\le n$. 
Using $f_{G_{j_0+1}} (z) \le 5+n\ep$ we obtain 
\begin{equation*}
\begin{split}
|g(z)| &\le \bigg| \sum_{j=1}^{j_0} f_{G_j} (z)\bigg| 
+ |f_{G_{j_0+1}} (z)| 
+\bigg| \sum_{j>j_0+1}^n f_{G_j} (z)\bigg|\\
&< 1+5+n\ep + (n-(j_0+1))n \ep < 6+n^2\ep\ .
\end{split}
\end{equation*}
Hence $|z| \le 7$ follows and the lemma is proved since 
$n/7>C$ and such vectors $y$ 
and $z$ can be produced in any subsequence of $(y_i)$.
\end{proof}



Our next lemma follows from Proposition 3.2 in \cite{AOST}.

\begin{lem}\label{lem5}
Let $X$ be a Banach space. 
Assume that for all $n$, $(x_i^n)_{i=1}^\infty$ is a normalized weakly null 
sequence in $X$ having spreading model $(\tilde x_i^n)$ which is not 
equivalent to the unit vector basis of $\ell_1$. 
Then there exists a normalized weakly null sequence $(y_i)\subseteq X$ 
with spreading model $(\tilde y_i)$ such that $(\tilde y_i)$ is not 
equivalent to the unit vector basis of $\ell_1$. 
Moreover, for all $n$ 
$$2^{-n} \| \sum a_i \tilde x_i^n\| \le \|\sum a_i \tilde y_i\|$$
for all $(a_i)\subseteq \real$.
\end{lem}

\begin{thm}\label{thm6}
If $X$ is elastic (and separable) then $c_0$ embeds into $X$.
\end{thm}

We postpone the proof to complete first  the 

\begin{proof}[Proof of  Theorem \ref{thm:elastic}]  
$\quad$\newline\indent
Assume that $X$ is infinite dimensional and every isomorph of $X$ is
$K$-elastic.  Then by Theorem~\ref{thm6}, $c_0$ embeds into $X$. 
Choose $k_n\uparrow\infty$ so that $2^{-n} k_n\to\infty$.
Using the renormings of $c_0$ by 
$$|(a_i)|_n = \sup \bigg\{ \Big|\sum_F a_i\Big| : |F| = k_n\bigg\}$$
and that $X$ is $K$-elastic we can find for all $n$ a normalized 
weakly null sequence $(x_i^n)_{i=1}^\infty \subseteq X$ with spreading 
model $(\tilde x_i^n)_{i=1}^\infty$ satisfying 
$$\|\sum a_i \tilde x_i^n\| \ge K^{-1} |(a_i)|_n\ $$
 and moreover each $(\tilde x_i^n)$ is equivalent to the unit vector basis 
of $c_0$. 
Thus by Lemma~\ref{lem5} there exists a normalized weakly null sequence 
$(y_i)$ in $X$ having spreading model $(\tilde y_i)$ which is not 
equivalent to the unit vector basis of $\ell_1$ and which satisfies 
for all $n$, 
$$\bigg\| \sum_1^{k_n} \tilde y_i\bigg\| \ge K^{-1} 2^{-n} k_n
\to \infty\ .$$
Thus $(\tilde y_i)$ is not equivalent to the unit vector basis of $c_0$ 
as well.

By Lemmas~\ref{lem1} and \ref{lem3}, for all $C<\infty$ we
can  find $n\in\nat$ and a renorming $Y$ of $X$ so that $Y$ contains a
normalized  weakly null sequence admitting no  subsequence which is block 
$n$-unconditional with constant $C$. By the assumption on $X$, the space
$Y$ must $K$-embed into every isomorph of $X$.  
But if $C$ is large enough this contradicts Lemma~\ref{lem2}.
\end{proof}

It remains to prove Theorem~\ref{thm6}. 
We shall employ an index argument involving $\ell_\infty$-trees 
defined on Banach spaces.
If $Y$ is a Banach space our trees $T$ on $Y$ will be countable.
For some $C$ the nodes of $T$ will be elements $(y_i)_1^n \subseteq Y$ 
with $(y_i)_1^n$ bimonotone basic and satisfying $1\le \|y_i\|$ and 
$\|\sum_1^n \pm y_i\| \le C$ for all choices of sign. 
Thus $(y_i)_1^n$ is $C$-equivalent  to the unit vector basis of 
$\ell_\infty^n$. 
$T$ is partially ordered by $(x_i)_1^n \le (y_i)_1^m$ if $n\le m$ and  
$x_i =y_i$ for $i\le n$. 
The order $o(T)$ is given as follows.  
If $T$ is not well founded (i.e., $T$ has an infinite branch), then 
$o(T) = w_1$. 
Otherwise we set for such a tree $S$, $S' = \{(x_i)_1^n\in S : (x_i)_1^n$ 
is not a maximal node$\}$.
Set $T_0 = T$, $T_1=T'$ and in general $T_{\alpha+1} = (T_\alpha)'$ 
and $T_\alpha = \cap_{\beta <\alpha} T_\beta$ if $\alpha$ is a limit ordinal. 
Then 
$$o(T) = \inf \{\alpha : T_\alpha = \phi\}\ .$$

By Bourgain's index theory \cite{B}, 
if $X$ is separable and contains for all $\beta < \w_1$  
such a tree of index at least $\beta$,  then $c_0$ embeds into $X$. 

We now complete the 

\begin{proof}[Proof of Theorem~\ref{thm6}]
$\quad$\newline\indent 
Without loss of generality we may assume that $X\subseteq Z$ where $Z$ 
has a bimonotone basis $(z_i)$.
Let $X$ be $K$-elastic.
We will often use semi-normalized sequences in $X$ which are a tiny 
perturbation of a block basis of $(z_i)$ and to simplify the estimates 
we will assume below that they are in fact a block basis of $(z_i)$. 

For example, if $(y_i)$ is a normalized basic sequence in $X$ then we call 
$(d_i)$ a {\em difference sequence\/} of $(y_i)$ if $d_i = y_{k(2i)} 
- y_{k(2i+1)}$ for some $k_1<k_2<\cdots$. 
We can always choose such a $(d_i)$ to be a semi-normalized perturbation 
of a block basis of $(z_i)$ by first passing to a subsequence $(y'_i)$ 
of $(y_i)$ so that $\lim_{i\to\infty} z_j^* (y'_i)$ exists for all $j$, 
where $(z_i^*)$ is biorthogonal to $(z_i)$, and taking $(d_i)$ to be a 
suitable difference sequence of $(y'_i)$. 
We will assume then that $(d_i)$ is in fact a block basis of $(z_i)$. 

We inductively construct for each {\em limit\/} ordinal $\beta < \w_1$, a 
Banach space $Y_\beta$ that embeds into $X$. 
$Y_\beta$ will have a normalized bimonotone basis $(y_i^\beta)$ that can 
be enumerated as $(y_i^\beta)_{i=1}^\infty = \{y_i^{\beta,\rho,n} 
: \rho\in C_\beta$,
$n\in \nat$, $i\in\nat\}$ where $C_\beta$ is some countable set. 
The order is such that $(y_i^{\beta,\rho,n})_{i=1}^\infty$ is a 
subsequence of $(y_i^\beta)$ for fixed $\rho$ and $n$.

Before stating the remaining properties of $(y_i^\beta)$ we need some 
terminology. 
We say that $(w_i)$ is a 
{\em compatible difference sequence of $(y_i^\beta)$ of order~1\/} 
if $(w_i)$ is a difference sequence of 
$(y_i^\beta)$ that can be enumerated as follows, 
$$(w_i) = \{ w_i^{\beta,\rho,n} : \rho\in C_\beta,\ n,i\in \nat\}$$
and such that for fixed $\rho$ and $n$, 
$$(w_i^{\beta,\rho,n})_i\ \text{ is a difference sequence of }\ 
(y_i^{\beta,\rho,n+1})_i\ .$$
If $(v_i)$ is a compatible difference sequence of $(w_i)$ of order~1, in the 
above sense, $(v_i)$ will be called a {\em compatible difference sequence 
of $(v_i)$ of order~2}, and  so on. 
$(y_i^\beta)$ will be said to have order~0.

Let $(v_i)$ be a compatible difference sequence of $(y_i^\beta)$ of some 
finite order. 
We set 
\begin{equation*}
\begin{split}
T\big((v_i)\big) & = \bigg\{ (u_i)_1^s : \text{ the $u_i$'s are distinct 
elements of $\{v_i\}_1^\infty$, possibly in different order,}\\
&\qquad \text{and }\ \bigg\|\sum_1^s \pm u_i\bigg\|=1
\text{ for all choices of sign}\bigg\}\ .
\end{split}
\end{equation*}
$T\big((v_i)\big)$ is then an $\ell_\infty$-tree as described  above with
$C=1$.  The inductive condition on $Y_\beta$, or should we say on 
$(y_i^{\beta,\rho,n})$, is that for all compatible difference sequences 
$(v_i)$ of $(y_i^{\beta,\rho,n})$ of finite order, 
$$o(T((v_i))) \ge \beta\ .$$
\renewcommand{\qed}{}
\end{proof}

Before proceeding we have an elementary

\begin{sublem}
Let $C<\infty$ and let $(w_i)$ be a block basis of a bimonotone basis $(z_i)$
with 
$1\le \|w_i\|\le C$ for all $i$ and let 
$$\A = \{F\subseteq \nat :F \text{is finite and }
\|\sum_{i\in F}\pm w_i\| \le C
\text{ for all choices of sign}\}.$$
Then there exists an equivalent norm $|\cdot|$ on $[(w_i)]$ so that 
$(w_i)$ is a bimonotone normalized basis such that for all $F\in\A$, 
$$\Big| \sum_F \pm w_i\Big| = 1\ .$$
\end{sublem}

\begin{proof}
Define $|\sum a_i w_i| = \|(a_i)\|_\infty \vee C^{-1}\|\sum a_i w_i\|$.
\end{proof}

We begin by constructing $Y_w$. 
Let $(x_i)\subseteq X$ be a normalized block basis of $(z_i)$. 
For $n\in\nat$,  let $|\cdot|_n$ be an equivalent norm on $[(x_i)]$ 
given by the sublemma for $C= 2^n$. 
Thus $|\sum_F \pm x_i|_n =1$ if $|F| \le 2^n$.

Since $X$ is $K$-elastic, 
for all $n$, $([(x_i)],|\cdot|_n)$ $K$-embeds into $X$. 
We thus obtain for $n\in\nat$, a sequence $(x_i^n)_i\subseteq X$ with 
$1\le \|x_i^n\| \le K$ for all $i$ and such that 
$1\le \|\sum_{i\in F} \pm x_i^n\| \le K$ for all $|F|\le 2^n$ and 
all choices of sign.
Furthermore $(x_i^n)_i$ is $K$-basic. 
By standard perturbation and diagonal arguments we may for each $n$
pass to a difference sequence $(d_i^n)_i$ of $(x_i^{n+1})_i$ so 
that enumerating, $(d_i) = \{d_i^n :n,i\in\nat\}$ is a block basis 
of $(z_i)$ with $1\le \|d_i\| \le K$ and 
with each $(d_i^n)_i$ being a subsequence of $(d_i)$. 
We have that for $|F|\le 2^n$ and all signs, 
$$1\le \bigg\| \sum_{i\in F}\pm d_i^n\bigg\| \le K\ .$$

We renorm $[(d_i)]$ by the sublemma for $C=K$ and let the ensuing  space 
be $Y_\w$. 
We change the name of $(d_i)$ to $(y_i^w)$ in this new norm and let 
$(y_i^{\w,1,n})_i = (d_i^n)_i$.
$(y_i^\w)$ has the property that if $(w_i)$ is a compatible difference 
sequence of $(y_i^\w)$ of finite order, then $o(T((w_i))) \ge \w$. 
Indeed, if $(w_i) = \{w_i^{\w,1,n} :n\in\nat$, $i\in\nat\}$ then for 
$|F|\le 2^n$, $\|\sum_{i\in F} \pm w_i^{\w,1,n}\| =1$. 

Assume that $Y_\beta$ has been constructed for the limit ordinal $\beta$
with basis $(y_i^\beta)= \{y_i^{\beta,\rho,n} :p\in C_\beta$, 
$n,i\in\nat\}$ with the requisite properties above.

Let $U: Y_\w\to X$ and $V:Y_\beta \to X$ be $K$-embeddings. 
Since in total we are dealing with a countable set of sequences, namely 
$(y_i^{\w,1,n})_i$ for $n\in\nat$ and $(y_i^{\beta,\rho,n})_i$ 
for $p\in C_\beta$, $n\in\nat$, by diagonalization and perturbation we 
can find a compatible difference sequence $(w_i^\w)_i$ of $(y_i^\w)_i$ of 
order~1 and a compatible difference sequence $(w_i^\beta)_i$ of 
$(y_i^\beta)_i$ of order~1 so that under a suitable reordering, 
$(d_i)= (U w_i^\w)_i \cup (Vw_i^\beta)_i$ is a block basis of $(z_i)$. 
Moreover each $(Uw_i^{\w,1,n})_i$ and $(Vw_i^{\beta,\rho,n})_i$ is a  
subsequence of $(d_i)_i$.

Adjoin a new point $p_0$ to $C_\beta$ and set $C_{\beta+\w} = C_\beta
\cup \{\rho_0\}$.  
Let $d_i^{\beta+\w,\rho,n} = Vw_i^{\beta,\rho,n}$ for $\rho\in C_\beta$ 
and $d_i^{\beta+\w,\rho_0,n} = Uw_i^{\w,1,n}$. 
For $|F| \le 2^n$ and $G\subseteq C_\beta\times \nat\times \nat$ for which  
$\big\|\sum\limits_{(\rho,n,i)\in G} \pm w_i^{\beta,\rho,n}\big\|=1$, 
we have by the triangle inequality that 
$$\bigg\|\sum_{i\in F} \pm d_i^{\beta+\w,\rho_0,n} 
+ \sum_{(\rho,n,i)\in G} \pm d_i^{\beta+\w,\rho,n} \bigg\| \le 2K.$$

It follows that if we let $(y_i^{\beta+\w})$ be the basis $(d_i^{\beta+\w})$, 
renormed by the sublemma for $C=2K$, that $Y_{\beta+\w} = [(y_i^{\beta+\w})]$
has the required properties.

If $\beta$ is a limit ordinal not of the form $\alpha+\w$ we let 
$U_\alpha : Y_\alpha \to X$ be a $K$-embedding for each limit ordinal 
$\alpha <\beta$.
We again diagonalize to form compatible difference sequences 
$(w_i^\alpha)_i$ of $(y_i^\alpha)$ of order~1 for each such $\alpha$ 
so that $(U_\alpha w_i^\alpha)_{\alpha,i}$ is a block basis of $(z_i)$ 
in some order. 
We let $C_\beta$ be a disjoint union of the $C_\alpha$'s and in the 
manner above obtain $Y_\beta$.\qed

$C[0,1]$ is, of course, $1$-elastic. 
By virtue of Lemma~\ref{lem2}, for all $K$, $C[0,1]$ can be renormed to 
be elastic but not $K$-elastic. 
Are there other examples of separable elastic spaces?

\begin{problem}\label{prob7}
Let $X$ be elastic (and separable, say). 
Does $C[0,1]$ embed into $X$?
\end{problem}

Using index arguments, we have the following partial result.

\begin{prop}\label{prop:elastic}
Let  $X$ be a separable Banach space, and 
suppose that $Y=\sum X_n$ is a symmetric decomposition of a space $Y$
into spaces uniformly isomorphic to $X$.  If
$Y$ is elastic, then $C[0,1]$ embeds into $X$.  

In particular, if $1\le p<\infty$ and $\ell_p(X)$ is elastic, then 
$C[0,1]$ embeds into $X$.
\end{prop}
\begin{proof}

 Let us first observe that if  $C[0,1]$ embeds
into $Y$, then $C[0,1]$ embeds into $X$. 
Since this is surely well known, we just sketch a proof (which, incidentally,
uses only that the decomposition $Y=\sum X_n$  of $Y$ is  unconditional): Let
$P_n$ be the projection from
$Y$ onto $X_n$ and let $W$ be a subspace of $Y$ which is isomorphic to
$C[0,1]$. By a theorem of Rosenthal's \cite{R}, it is enough to show that for
some $n$, the adjoint ${P_n}_{|W}^*$  of the restriction of $P_n$ to $W$ has
non separable range.  This will be true if there is an $m$ so that
${S_m}_{|W}^*$ has non separable range, where $S_m:=\sum_{i=1}^m P_i$. Let $Z$
be a subspace of $W$ which is isomorphic to $\ell_1$. If no such $m$ exists,
then for every
$m$, the restriction of $S_m$ to $Z$ is strictly singular (that is, not an
isomorphism on any infinite dimensional subspace of $Z$), and it then follows
that
$Z$ contains a sequence
$(x_n)$ of unit vectors which is an arbitrarily small perturbation of a
sequence
$(y_n)$ which is disjointly supported.  The sequence $(y_n)$, a fortiori 
$(x_n)$, is then easily seen to be equivalent to the unit vector basis of
$\ell_1$ and its closed span is complemented in $Y$ since the decomposition is
unconditional.  It follows that
$\ell_1$ is isomorphic to a complemented subspace of $C[0,1]$, which of course
is false.

To complete the proof of Proposition~\ref{prop:elastic}, we assume that 
$Y$ is $K$-elastic and prove that $C[0,1]$ embeds into $Y$.  The proof is
similar to, but simpler than, the proof of Theorem~\ref{thm6}.  
First we recall the definition of certain 
canonical trees $T_\alpha$ of order $\alpha$ for $\alpha<\omega_1$
(see e.g. \cite{JO}). 
These form the frames upon which we will hang our bases.  
The tree $T_1$ is a single node.  
If $T_\alpha$ has been defined, we choose a new node
$z\not\in T_\alpha$ and set $T_{\alpha + 1}:= T_{\alpha}\cup\{z\}$, ordered by 
$z< t$ for all $t\in T_\alpha$ and with $T_\alpha$ preserving its order.  
If $\beta<\omega_1$ is a limit order, we let $T_\beta$ be the disjoint 
union of $\{T_\alpha : \alpha<\beta\}$.  
Then if $s,t$ are in $T_\beta$, we say that 
$s\le t$ if and only if $s,t$ are both in $T_\alpha$ for some $\alpha<\beta$
and $s\le t$ in $T_\alpha$.  

We shall prove by transfinite induction that if $(x_i)_{i=1}^\infty$ is any
  normalized monotone basic sequence and $\alpha<\omega_1$, there there is a
Banach space
$Y_\alpha=Y_\alpha(x_i)$ with a  normalized monotone basis
$(y_t^\alpha)_{t\in T_\alpha}$ so that  if $(\gamma_i)_{i=1}^n$ is any branch
in
$T_\alpha$ then $(y_{\gamma_i}^\alpha)_{i=1}^n$ is $1$-equivalent to
$(x_i)_{i=1}^n$.  Furthermore, each  $Y_\alpha$ will $K$-embed  into $Y$. Just
as in the proof of Theorem~\ref{thm6}, it then follows from index theory that
the Banach space spanned by  $(x_i)_{i=1}^\infty$ embeds into
$Y$.  

Fix any  normalized monotone basic sequence $(x_i)_{i=1}^\infty$. 
Suppose that $\beta$ is a limit ordinal and
$Y_\alpha=Y_\alpha(x_i)_{i=1}^\infty$ has been defined for all
$\alpha<\beta$.   In view of the hypotheses, the space $Y$ has a symmetric
decomposition into spaces uniformly isomorphic $Y$, which we can index as
$Y=\sum\limits_{\alpha<\beta} X_\alpha$. For each $\alpha$, there is an
isomorphism $L_\alpha$ from $Y_\alpha$ into $X_\alpha$ so that $\|L_\alpha\|=1$
and $\|L_\alpha^{-1}\|$ is bounded independently of $\alpha$.  We can put an
equivalent norm on 
$Y=\sum\limits_{\alpha<\beta} X_\alpha$ to make each $L_\alpha$ an isometry
and make the decomposition $1$-unconditional (but not necessarily
$1$-symmetric).  Define $Y_\beta$ to be the closed linear span of 
$\{L_\alpha Y_\alpha : \alpha<\beta\}$ in $Y$ with its new norm.  The space
$Y_\beta$ has the desired basis indexed by $T_\beta$ and $Y_\alpha$ must
$K$-embed into $Y$ with its original norm because $Y$ is $K$-elastic.

If $\beta=\alpha + 1$, we let 
$Y_\beta(x_i)_{i=1}^\infty=\real \oplus Y_\alpha (x_i)_{i=2}^\infty$
with the norm given by
$$
\|(a,y)\|:= 
\sup\{\|a x_1+\sum_{i=2}^n y_{\gamma_i}^{\alpha *}(y) y_{\gamma_i}^{\alpha}\|,
$$
where the supremum is taken over all $(\gamma_i)_{i=2}^n$ which form a branch
or an initial segment of a branch in $T_\alpha$.

Again it is clear that $Y_\beta$ must $K$-embed into $Y$ and that $Y_\beta$
has the desired basis.  (In the case $\beta=\alpha+1$, the space $Y_\beta$ need
not contain
$Y_\alpha$ isometrically, but that is irrelevant.) 
\end{proof}

\end{section}

\end{document}